\documentclass[a4paper,10pt,leqno]{article}
\usepackage[applemac]{inputenc}
\usepackage[T1]{fontenc}
\usepackage[english,francais]{babel}

\usepackage{geometry}
\usepackage{amsmath}
\usepackage{amssymb,amsmath,amsthm,amscd}
\usepackage{mathrsfs}

\usepackage{sectsty}
\usepackage{lipsum}

\allsectionsfont{\centering}
\usepackage{enumerate}

\theoremstyle{plain} 
\newtheorem{theo}{Théorème}[section]

\newtheorem{defi}{Définition}
\newtheorem{rema}{Remarque}
\newtheorem{exem}{Exemple}
\newtheorem{prop}{Proposition}[section]

\newtheorem{hypo}{Hypothèse}

\numberwithin{equation}{section}

\title{L'INÉGALITÉ DE MELIN-HÖRMANDER EN CARACTÉRISTIQUES MULTIPLES}
\author{BERNARD LASCAR AND RICHARD LASCAR}
\date{}
\newcommand{\R}{\mathbb{R}}

\newcommand{\N}{\mathbb{N}}
\newcommand{\LL}{\mathcal{L}}
  
\begin{document}
\maketitle
\begin{center}
\it Dédié à la mémoire du Professeur Louis Boutet de Monvel
\end{center}
{\abstract On prouve ici l'inégalité de Melin-Hörmander en caractéristiques $2k$ pour $k\in\N$ quelconque.

We prove in this work the Melin-Hörmander inequality for operators with multiple characteristics.

AMS Classification : 35S05.}
\section{NOTATIONS}
On considère dans ce travail des opérateurs pseudo-différentiels classiques sur $\R^n$, $P=P(x,D)$ où :
\begin{equation}\label{eq:defpux} Pu(x)=\int e^{ix\xi} p(x,\xi)\widehat{u}(\xi)\frac{d\xi}{(2\pi)^n},
\end{equation}
où $p(x,\xi)$ est classique, c'est à dire admet un développement :
\begin{equation}
p(x,\xi)=p_m(x,\xi) + \ldots + p_{m-j}+\ldots
\end{equation}
où $p_{m-j}$ est positivement homogène en $\xi$ de degré $m-j$ ; $m\in\R$ est l'ordre de l'opérateur ; et $u\in\mathcal S(\R^n)$.

Dans la mesure où on va produire une inégalité de borne inférieure, on préférera à la formule (\ref{eq:defpux}) la formule de la quantification de Weyl :
\begin{equation}\label{eq:quanWey} Pu(x)=(\sigma^w) u(x)=\int e^{i(x-y)\xi} \sigma\Bigl(\frac{x+y}{2},\xi\Bigr) u(y)dy\frac{d\xi}{(2\pi)^n}.
\end{equation}
Comme :
\begin{equation} \sigma(x,\xi)=e^{-i(D_x,D_{\xi})/2}p(x,\xi);\;\sigma\text{ est aussi un symbole classique de degré }m.
\end{equation}
Par ailleurs $((\sigma)^w)^*=(\overline{\sigma})^w$ et $\sigma_{m-1}$ est le symbole sous principale de $P$ bien défini sur $\Sigma=\{(x,\xi)\in\R^n\times \R^n; p_m(x,\xi)=dp_m(x,\xi)=0\}$.

Bien que notre résultat ne concerne que des symboles classiques nous aurons besoin des résultats du calcul de Weyl d'Hörmander [7] volume III page 150. Ce calcul concerne des symboles très généraux dont les propriétés sont définies par un pois $m(X)$ et une métrique $g_X$ où maintenant $X$ désigne le point courant de $\R^n\oplus \R^n$.
\begin{defi} Soit $X\in\R^n\to g_X$ un champ de métriques, on dit que $g_X$ est admissible si on a les trois conditions :

lenteur 
\begin{equation} \text{Il existe }c \text{ et } C \text{ positives telles que }g_X(Y)\leq c \text{ entraîne } g_{X+Y} \leq Cg_X.
\end{equation}
tempérence 
\begin{equation} \text{Il existe }c \text{ et } N \text{ tels que }g_X\leq Cg_Y(1+g_X^\sigma(X-Y))^N.
\end{equation}
principe d'incertitude 
\begin{equation} h^2(X)=\sup_{v\in\R^n\oplus\R^n} \frac{g_X(v)}{g_X^\sigma(v)}\leq 1.
\end{equation}
\end{defi}
Un pois $m(X)$ est admissible s'il est lent et admissible.

On dira que $P$ est dans la classe $S(m,g)$ si $p\in S(m,g)$ c'est à dire :
\begin{equation} |D_x^j p(X)(t_1,\ldots, t_j)|\leq C_j m(X) (g_X(t_1)\ldots g_X(t_j))^{1/2} \quad \text{pour tout } t_1,\ldots,t_j\in \R^n\oplus \R^n.
\end{equation}
Les opérateurs se composent et $p_1\circ p_2=(p_1 \# p_2)^w$, $p_1\# p_2=q$ avec :
\begin{equation}
q(X)-\sum_{j<N} (\tfrac{1}{2})^j \frac{1}{j !} \sigma (D_x,D_\xi ; D_y, D_\eta)^j (p_1(x,\xi)p_2(y,\eta))\vert_{(x,\xi)=(y,\eta)}\in S(m_1m_2 h^N, g).
\end{equation}
Et $q\in S(m_1m_2, g)$ si $p_1\in S(m_1,g)$ et $p_2\in S(m_2,g)$. On a : 
\begin{theo} Si $a\in S(1,g)$ alors :
\begin{equation} \|a^w\|_{(L(\R^n,\R^n))}\leq C_a \text{ où } C_a \text{ est une semi-norme de } S(1,g).
\end{equation}
\end{theo}
\begin{theo} L'inégalité de Gärding (sharp) s'énonce : si $a\in S(h^{-1},g;\LL(F,F))$ où $F$ est un espace de Hilbert et si pour tout $X\in \R\oplus \R^n$, $a(X)\geq 0$ alors il existe $C_a$ ne dépendant que des semi-normes de $a$ tel que 
\begin{equation} (a+C_a I d)^w\geq 0.
\end{equation}
\end{theo}
\begin{theo} L'inégalité de Fefferman-Phong s'énonce : si $a\in S(h^{-2},g)$ est scalaire et si pour tout $X\in\R^n \oplus \R^n$, $a(X)\geq 0$ alors il existe $C_a$ ne dépendant que des semi-normes de $a$ tel que
\begin{equation} (a+C_a)^w\geq 0.
\end{equation}
\end{theo}
L'inégalité de Fefferman-Phong est très intéressante dans les cas où il y a pas du tout de calcul symbolique asymptotique (elle a d'ailleurs été inventée pour $(\psi)$) (voir [5] et [11]) mais n'est pas optimale en effet :
\begin{equation} a=x^2+\xi^2-1 \not\ge 0 \text{ et pourtant } a^w =D_x^2+x^2-1\geq 0
\end{equation}
En effet la fonction $h$ de ce calcul symbolique est $h(X)=\langle x,\xi\rangle^{-2}$, donc seules des perturbations négatives de degré $\langle x,\xi\rangle^{-2}$ ne seraient admises.

A. Melin puis L. Hörmander ont répondu à cette question, à l'aide des invariants symplectiques issus des travaux de L. Boutet de Monvel [3] et V. Ivrii-V M. Petkov [9] qui sont le symbole sous principal et la trace plus du Hessien.

En effet si $p_m\geq 0$, sur $\Sigma=\{(x,\xi)\in T^*\R^n\setminus 0 ; p_m(x,\xi)=0\}$, en plus de $p_{m-1}^s$ on dispose de la matrice fondamentale définie par la relation :
\begin{equation} \nabla^2 p_m(\rho)(t,t') =\sigma(t,F_{p_m}(\rho)t')\text{ pour } t \text{ et } t'\in T_{\rho}(T^*\R^n)
\end{equation}
On voit assez facilement que $\text{Spec } F_{p_m}(\rho)\subset i \R$ et est formé de couples de valeurs propres $\pm i\lambda_j,\lambda_j\geq 0$. Alors on note :
\begin{equation} \text{tr}^+ P_{p_m}(\rho)=\sum_j\lambda_j.
\end{equation}
L. Hörmander a prouvé dans [6] le théorème suivant :
\begin{theo}\label{theo:Horm1}
Si $P(x,D_x)$ est un opérateur pseudo-différentiel classique auto-adjoint sur $\R^n$ qui satisfait à :
\begin{hypo} \begin{enumerate}[i)] \item $p_{2m}\geq 0$ et $p_{2m}^{-1}(0)=\Sigma$, où $\Sigma$ est une variété lisse de $T^* \R^n\setminus 0$ de rang symplectique constant. On note $d_\Sigma$ la distance à $\Sigma$.
\item 
\begin{equation} p_{2m}\simeq d_\Sigma^2 \text{ pour } (x,\xi)\in S^*(\R^n).
\end{equation}
\item En chaque point $\rho\in\Sigma$,
\begin{equation}\label{eq:p2m-1s1} p_{2m-1}^s(\rho)+\frac{1}{2}\text{tr}^+ F_{p_{2m}}(\rho)\geq 0.
\end{equation}
Alors on a :
\begin{equation}\label{eq:Puu}
(Pu,u)\geq -C_K\|u\|_{(m-1)}^2 \text{ pour }u\in C_K^{\infty}
\end{equation}
\end{enumerate}
\end{hypo}
\end{theo}
Bien entendu si (\ref{eq:p2m-1s1}) est remplacée par :
\begin{equation}\label{eq:p2m-1s2}
p_{2m-1}^s(\rho) +\frac{1}{2}\text{tr}^+ F_{p_{2m}}(\rho) >0,
\end{equation}
alors (\ref{eq:Puu}) devient 
\begin{equation}\label{eq:Puugeq}
(Pu,u)\geq c_K\|u\|_{(m-1)/2}^2 -C_K\|u\|_{m-1}^2 \text{ pour } u\in C_K^{\infty} ;
\end{equation}
qui s'avère en pratique plus utile.
\section{ÉNONCÉ DU RÉSULTAT}
On étend à la multiplicité d'ordre supérieure le théorème  (\ref{theo:Horm1}) ou plus exactement sa conséquence (\ref{eq:Puugeq}) quand on a la condition (\ref{eq:p2m-1s2}). Il est clair que l'invariance exige lorsque le symbole principal s'annule à l'ordre $2k$ sur une variété $C^{\infty}\Sigma$ que pour $0\leq j\leq k$, $p_{m-j}$ s'annule à l'ordre $2k-2j$ sur $\Sigma$. On peut alors énonce notre résultat.
\begin{hypo}\label{hypo:2} Soir $P(x,D_x)$ un opérateur pseudo-différentiel classique quto-ajkoint de degré $m$ sur $\R^n$.
\begin{enumerate}[i)] \item $p_m\geq 0$ et $p_m^{-1}(0)=\Sigma$, où $\Sigma$ est une variété lisse de $T^* \R^n\setminus 0$ de rang symplectique constant. On note $d_\Sigma$ la distance à $\Sigma$.
\item Il existe $k\in\N$, tel que :
\begin{equation}\label{eq:2.1} p_{m-j}=\mathcal O(d_\Sigma^{2k-2j}) \text{ pour }O\leq j\leq k \text{ et } p_m\simeq d_\Sigma^{2k} \text{ pour }(x,\xi)\in ST^*(\R^n).
\end{equation}
\item En chaque point $p\in\Sigma$, l'opérateur différentiel localisé $P_\Sigma$ a une borne inférieure strictement positive sur $L^2$.
\end{enumerate}
\end{hypo}
Les hypothèses \ref{hypo:2} sont invariantes par transformation canonique homogène.

On déduit des hypothèses (\ref{hypo:2}) l'inégalité :
\begin{theo}\label{theo:2.1} Soit $P(x,D_x)$ un opérateur pseudo-différentiel classique auto-adjoint de degré $m$ qui vérifie les hypothèses (\ref{hypo:2}), alors pour tout compact $K$ de $\R^n$, il existe des constantes $c_K$ et $C_K$ positives telles que :
\begin{equation} \label{eq:2.2}
(Pu,u) \geq c_K\|u\|_{m/2-k/2}^2 -C_K\|u\|_{m/2-k/2-1/2}^2\text{ pour tout }u\in C_0^\infty(K).
\end{equation}
\end{theo}
Le théorème (\ref{theo:2.1}) est exactement l'extension à $k$ quelconque du théorème 22.33 page 364, volume III du traité de L. Hörmander [7] qui traite le cas $k=1$.

L'opérateur localisé a été introduit par L. Boutet de Monvel [3], les invariances par transformations canoniques ont été largement développées à cette époque.

On donne une preuve moderne qui s'appuie sur un argument de deuxième microlocalisation.

Nous allons d'abord définir correctement $P_\Sigma$. Pour cela on reprend l'exposé de Boutet de Monvel-Grigis-Helffer [4]. On travaille dans des coordonnées $(u,v)$ près de $\rho\in\Sigma$ où $\Sigma=\{u=0\}$. L'invariance par changement de coordonnées est une autre question.
\begin{defi} On dit que $P\in\mathcal N^{m,k}$, si $P(x,D)$ est une opérateur pseudo-différentiel classique de degré $m$ tel que pour $0\leq j\leq k/2$, $p_{m-j}$ s'annule à l'ordre $k-2j$ sur $\Sigma=\{u_j=1,1\leq j\leq p\}$ et $T_k(P)$ désigne la somme obtenue en prenant les développements de Taylor sur $u_j=0$, ainsi obtenus.
\end{defi}
On décompose les champs de vecteurs $\frac{\partial}{\partial x_S}$ et $\frac{\partial}{\partial \xi_S}$ suivant leurs composantes tangentes et normales à $\Sigma$ :
\begin{equation} \frac{\partial}{\partial x_S}=\sum_{j} B_{js}(v)\frac{\partial}{\partial u_j}+r_S
\end{equation}
\begin{equation} \frac{\partial}{\partial \xi_S}=\sum_{j} C_{js}(v)\frac{\partial}{\partial u_j}+\rho_S
\end{equation}
où $r_S$ et $\rho_S$ sont des champs tangents à $\Sigma$.
L'opérateur $B(v)$ de matrice $B_{j,s}$ opère de $E=\R^n$ dans $N=\R^p$; $C : E^n\to N$ de matrice $C_{js}$. $L=(B+C):E\times E^* \to N$, puis enfin $A=C^t B : N^*\to N$.
\begin{defi} Si $L(x,\xi)=Bx+C\xi$, et $A=C^tB\in\LL(N^*,N)$ on pose
\begin{enumerate}[i)]
\item
\begin{equation} a_Lf=\int e^{ix\xi}a(L(x,\xi))\widehat f(\xi) \frac{d\xi}{(2\pi)^n}\text{ pour }a\in \mathcal S(N).
\end{equation}
\item On pose :
\begin{equation} p_\Sigma =a_A\text{ avec } a=T_kp
\end{equation}
\end{enumerate}
\end{defi}
On a dans [4] la propriété suivante :
\begin{prop} 
\begin{equation}
(P\circ Q)_\Sigma=P_\Sigma\circ Q_\Sigma
\end{equation}
\end{prop}
Et donc pour calculer $P_\Sigma$ il faut décomposer $P$ en produit d'opérateurs :
\begin{equation}
P=\sum_{|\alpha|\leq k}A_\alpha U_1^{\alpha_1}\ldots U_p^{\alpha_p}
\end{equation}
Si $P=U_1$ si $u_1=x_1$, $P_\Sigma=x_1$ ; Si $u_1=\xi_1$, $P_\Sigma=\frac{\partial}{i\partial_{x_1}}$.

Maintenant on veut faire le lien avec la quantification de Weyl.

On a besoin déjà d'un calcul symbolique, on quantifie donc les symboles avec un grand paramètre $\Lambda$ par la formule :
\begin{equation}
q^{w_\Lambda}u(x)=\int e^{i\Lambda(x-y)\xi} q\Bigl(\tfrac{x+y}{2},\xi,\Lambda\Bigr)u(y)dy\frac{\Lambda^n d\xi}{(2\pi)^n}
\end{equation}
On note pour $I$ un $k$-uple $I=(i_1,\ldots,i_k)$, $U^{(I)}$ le symétrisé pour toutes les permutations des indices du produit $U_{i_1}\circ \ldots \circ U_{i_k}$.
Soit $m\in\R$ et $\alpha\in\N$, définissons alors $T^{m,\alpha}$ comme la classe des fonctions qui s'écrivent :
\begin{equation} 
a=\sum_{2\leq 2l\leq \alpha} a_l u^{-2l +\alpha} \text{ où } a_l\in S(\Lambda^{-l},\Gamma), \Gamma=|dX|^2.
\end{equation}
Alors si $a\in T^{m,l}$ et $M\in\N$, $\Lambda^{M/2}\nabla^M a\in T^{m,l}$, donc $T^{m,l}\# T^{m',l'}\subset T^{m+m',l+l'}$ et $u^\alpha\in T^{0,\alpha}$.
\begin{prop} 
\begin{equation}\label{eq:2.11}
U^{(I)}=(u^{(I)})^w +r_I^{w} \text{ avec } r_I\in T^{0,|I|}.
\end{equation}
\end{prop}
La proposition se prouve aisément en utulisant la forme explicite du calcul de Weyl qui annule les temers impairs du développement lorsqu'on fait un produit symétrique. De plus si les $u_j$ sont des coordonnées le calcul (\ref{eq:2.11}) est exact et $r_I=0$.
\begin{exem} Supposons et ce n'est pas restrictif que $P\in\mathcal N^{0,2}$ soit écrit de façon symétrique :
\begin{equation} 
P=\alpha x^2 +\beta\Bigl( x\frac{D_x}{\Lambda} +\frac{D_x}{\Lambda}x\Bigr) +\gamma\Bigl(\frac{D_x}{\Lambda}\Bigr)^2 +\Lambda^{-1}p_1
\end{equation}
La condition $P_\Sigma$ signifie exactement (\ref{eq:2.1}) que :
\begin{equation}
Q=(\alpha_0 x^2 +2\beta_0 x\xi +\gamma_0 \xi^2 +\Lambda^{-1} p_1(0))^{w_\Lambda}>0\text{ sur }L^2(\R)
\end{equation}
Ce qui dans le cas des formes quadratiques se calcule sous la forme :
\begin{equation}
p_1(0) +\frac{1}{2} \text{tr}^{+} F_{Q_0} >0.
\end{equation}
\end{exem}
Calculé en coordonnées [4], on a :
\begin{equation}
P_\Sigma=\sum_{|\alpha|+|\beta|+2j=k}\frac{1}{\alpha !\beta !}\partial_x^{\alpha} \partial_\xi^{\beta} \sigma_{m-j}(x,\xi)(y^\alpha \eta^\beta)^w.
\end{equation}
\section{PREUVE DU THÉORÈME (\ref{theo:2.1})}
On ne va donner une preuve que dans le cas où $\Sigma$ est symplectique, les autres cas se déduisent directement. D'ailleurs les termes involutifs n'améliorent pas les conditions sur les termes d'ordres inférieurs et disposent d'un calcul symbolique avec $h$ petit donc ajouter des termes involutifs est une généralisation facile.

On commence par des réductions standards.

Le théorème (\ref{theo:2.1}) étant asymptotique et microlocal, il se réduit à une inégalité avec grand paramètre $\Lambda$ pour un opérateur pseudo-différentiel $P=(p)^{w_\Lambda}$ dont le symbole est classique :
\begin{equation} p(x,\xi,\Lambda)=\Lambda^{m} q_m(x,\xi)+\ldots +\Lambda^{m-1}q_{m-1}(x,\xi) +\ldots
\end{equation}
les $q_{m-j}$ peuvent aussi dépendre de $\Lambda$ avec $q_{m-j}\in S(1,\Gamma)$ où $\Gamma=|dX|^2$ et $X=(x,\xi)$.

L'inégalité \ref{eq:2.2} est écrite microlocalement :
\begin{prop} Il existe un petit voisinage $V(\rho_0)$ de chaque point $\rho_0$ de $\Sigma$ et pour chaque $N\in\N$, des constantes $c>0$ et $C_N>0$ telles que :
\begin{equation}\label{eq:3.2} (P\Theta u,\Theta u)\geq c\Lambda^{m-k}\|\Theta u\|^2 -C_N\Lambda^{-N}\|u\|^2,
\end{equation}
où $\Theta=(\theta)^{w_\Lambda}.\theta(X)\in C_0^{\infty}(V(\rho_0))$.
\end{prop}
Vu les hypothèses sur $\Sigma$, on peut choisir des coordonnées symplectiques dans lesquelles $\Sigma=\{(x,\xi)\in \R^{2n}; x_1=\ldots = x_d=\xi_1=\ldots =\xi_d =0\}$. Puisque toutes les hypothèses sont invariantes par transformation canonique homogène.

On change les notations et on écrit maintenant $\widetilde X=(X,X_n)$ où $\widetilde X$ sont touts les variables de $\R^n\oplus \R^n$, $X\in\R^d\oplus\R^d$, $X_n\in\R^{n-d}\oplus \R^{n-d}$, $\Sigma=\{\widetilde X;X=0\}$.

L'opérateur localisé $P_\Sigma$ est alors un opérateur différentiel à symbole polynomial en $X$ de degré au plus $2k$. Quantifié en $w_\Lambda$, le symbole total n'a pas la bonne homogénéité, puisque $x$ et $\xi$ sont respectivement les symboles de $x$ et $\frac{D_x}{\Lambda}$, leur commutateur vaut $\Lambda^{-1}$ tout comme $\Lambda^{-1}q_{m-1}$. On se simplifie la tâche si on quantifie autrement.

On introduit la transformation unitaire dans $L^2$
\begin{equation}\label{eq:3.3}
T_\Lambda u(x,x_n)=u(\Lambda^{-1/2}x,x_n)\Lambda^{-d/2}.
\end{equation}
\begin{equation}\label{eq:3.4} 
T_\Lambda q^{w_\Lambda} T_\Lambda^{-1}=(q(\Lambda^{-1/2}x,\Lambda^{1/2}\xi,x_n,\xi_n))^{w_\Lambda}=(q(\Lambda^{-1/2}X,X_n))^{w_{1,\Lambda}}.
\end{equation}
Dans la formule (\ref{eq:3.3}), $w_{1,\Lambda}$ est la quantification $1$ dans la variable $x$ et $\Lambda$ dans la variable $x_n$.

Si $Q_l\in S(1,\Gamma)$ et $Q_l=\mathcal O(d_\Sigma^l)$, alors par la formule de Taylor à l'ordre $l$, $Q_l=X^l Q_{1,l}(X,X_n)$ où $Q_{1,l}\in S(1,\Gamma)$.

Si maintenant on remplace $q(X,X_n)$ par $q(\Lambda^{-1/2}X,X_n)$ (\ref{eq:3.4}) ; $Q_l=X^l Q_{1,l}(X,X_n)$ devient $F_f(X,X_n)=\Lambda^{-l/2} Q_{1,l}(\Lambda^{-1/2}X,X_n)X^l$.
\begin{equation}\label{eq:3.5}
F_{1,l}=Q_{1,l}(\Lambda^{-1/2} X,X_n)\in S(1,\Gamma_\Lambda) \text{ avec } \Gamma_\Lambda =\Lambda^{-1}|dX|^2+|dX_n|^2 ;
\end{equation}
\begin{rema}\label{rema:1} Cependant maintenant $F_{1,l}$ est seulement supporté par une boule $B_R=\{X,|X|< R\Lambda^{1/2}\}$ où $R$ est petit et désigne le diamètre de $V(\rho_0)$.
\end{rema}
Ce que l'on veut montrer est que l'on peut localiser dans l'espace des $X$ de façon à réduire $F_{l}(X,X)$ à $F_{l,0}(X,X_n)=\Lambda^{-1/2} Q_{1,l}(0,X_n)X^l$
Soi si $a\geq 0$,
\begin{equation}
g_a=\frac{|dX|^2}{d_a(X)^2}+|dX_n|^2
\end{equation}
avec $d_a(X)=|X|+a$, $X^l\in S(d_a^l,g_a)$.

La condition $S(1,\Gamma_\Lambda)\subset S(1,g_a)$ entraîne que $d_a(X)\leq C\Lambda^{1/2}$, il faut donc avoir $|X|\leq C\Lambda^{1/2}$ et prendre $a\leq C\Lambda^{1/2}$, la première de ces deux conditions sera assurée si on a pris des symboles, dans le calcul de départ, à support près de $\rho_0\in\Sigma$. On veut de plus avoir un calcul pseudo-différentiel, il faut donc vérifier les conditions du calcul de Weyl :
\begin{prop} Si $a\geq 1$, la métrique $g_a$ est lente et $\sigma$ tempérée et les poids $d_a^m$ sont lents et $\sigma$ tempérés pour tout $m\in\R$. La fonction $h$ du calcul $(1,\Gamma)$ corespondant vaut : $h_a(X)=\max(d_a(X)^{-2},\Lambda^{-1})\leq 1$. On veut de plus composer une classe $S(m,g_a)$ avec une classe $S(m',g_b)$, il vaut vérifier que $h_{a,b}(X)=\max((d_a,d_b)^{-1},\Lambda^{-1})\leq 1$, soit $\alpha\#\beta\in S(mm',g_{\min(a,b)})$, si $\alpha\in S(m,g_a)$ et $\beta\in S(m',g_b)$.
\end{prop}
C'est en effet bien connu et parfaitement évident. On a donc montré :
\begin{prop} Il résulte des hypothèses d'annulation, $\widetilde F_{2k-2j}=\Lambda^{m-j}F_{2k-2j}\in S(\Lambda^{m-k} d_a^{2k-2j},g_a)$ dans $\{|X|\leq R\Lambda^{1/2}\}$ et si $1\leq a\leq C\Lambda^{1/2}$. On se trouve bien localisé dans une zone $\{|X|\leq R\Lambda^{1/2}\}$ à cause de la remarque (\ref{rema:1}).
\end{prop}
Il est clair qu'il va falloir mieux localiser que dans des boules de rayon $R\Lambda^{1/2}$, même si $R$ est petit.

Soit $\rho\geq 1$ à choisir, comme la métrique $g_\rho$ est lente, on peut faire des partitions de l'unité avec des symboles $S(1,g_\rho)$. Soit :
\begin{equation}
\chi_1^2(X)+\chi_2(X)^2=1,
\end{equation}
avec $\chi_j\in S(1,g_\rho)$, $\text{supp } \chi_1\subset B(0,\rho)$ et $\text{supp } \chi_2\subset B(0,\rho/2)^c$.

\noindent {\bf 3.1. La zone $X_2=\{X,|X|\geq 1/3\rho\}\cap \{|X|\leq R\Lambda^{1/2}\}$.}

On commence par minorer la contribution de la zone $X_2$.

La partie principale obtenue pour $j=0$ vérifie $\widetilde F_m \geq c\Lambda^{m-k} (\rho+|X|)^{2k}$, à cause de l'éllipticité transverse. Or pour $j>1$, $\widetilde F_{2k-2j}\leq C(|X|+\rho)^{2k-2j}$. L'opérateur total est elliptique positif dans la classe $S(\Lambda^{m-k}d_\rho^{2k},g_\rho)$ si $\rho$ est assez grand.

L'application de l'inégalité de Fefferman-Phong permet de voir d'emblée que si $Q=T_\Lambda P T_\Lambda^{-1}$,
\begin{equation}
(Q\chi_2 \chi_0 u,\chi_2\chi_0 u)\geq c\Lambda^{m-k} \|\chi_2\chi_0 d_\rho^k u\|^2-\Lambda^{m-k} C\max(\rho^{-2},\Lambda^{-1})^2\|u\|^2.
\end{equation}
Prenant les commutateurs avec $\chi_2$ on trouve que :
\begin{multline}\label{eq:3.9}
\Lambda^{-m+k}((Q\chi_2^2)^{w_{1,\Lambda}}\chi_0 u,\chi_0 u)\geq c\|\chi_2 d_\rho^k \chi_0 u\|^2- C\|d_\rho^{k-2}\chi_0 u\|^2-\\
C_N\max (\rho^{-1},\Lambda^{-1})^N\|u\|^2\text{ pour tout } N\in \N,
\end{multline}
où $\chi_0$ est une microlocalisation au sens usuel dans $V(\rho_0)$, mais comme (\ref{eq:3.9}) est valable avec n'importe quelle fonction $\chi_2$ supportée par $|X|\geq c\rho$, on peut remplacer le $d_\rho^{k-2}$ par $d_\rho^{k-1}$, car dans $X_2$ $d_\rho\sim d_1$, ceci modulo un $d_\rho^{-\infty}$ donc :
\begin{multline}\label{eq:3.10}
\Lambda^{-m+k} ((Q\chi_2^2)^{w_{1,\Lambda}}\chi_0 u,\chi_0 u)\geq c\|\chi_2 d_1^k \chi_0 u\|^2- C\|d_1^{k-2}\chi_0 u\|^2-\\
C_N\max (\rho^{-2},\Lambda^{-1})^N\|u\|^2\text{ pour tout } N\in \N.
\end{multline}
\noindent {\bf 3.2. La zone $X_1=\{X,|X|\leq 2\rho\}\cap \{|X|\leq R\Lambda^{1/2}\}$.} Il faut d'abord réduire chacun des symboles $p_{m-j}$ à son développement de Taylor à l'ordre $2k-2j$ sur $X=0$, puisque c'est là qu'est exprimée la condition sur le localisé. On écrit la formule de Taylor de $p_{m-j}$ à l'ordre $2k-2j$ avec reste intégral :
\begin{multline}\label{eq:3.11}
\Lambda^{m-j}p_{m-j}(\Lambda^{-1/2}X,X_n)=\frac{1}{(2k-2j)!}\Lambda^{m-k}(\nabla{2k-2j}p_{m-j})(0,X_n)X^{2k-2j}+\\ \int_0^1\frac{(1-\tau)^{2k-2j}}{(2k-2j)!}\Lambda^{m-k-1/2}(\nabla^{2k-2j-1}p_{m-j})(\tau \Lambda^{-1/2}X,X_n)X^{2k-2j+1}
\end{multline}
Dans le second membre de (\ref{eq:3.11}), aucun des deux termes n'a de support localisé en $X$, le premier terme est bien dans la formule du localisé, pour estimer correctement le second terme et ses dérivées on doit multiplier la formule (\ref{eq:3.11}) par une fonction de la première localisation notée $\chi_3$ qui vaut identiquement $1$ sur le support des $p_{m-j}$. Ceci s'écrit :
\begin{multline}\label{eq:3.12}
\Lambda^{m-j}p_{m-j}(\Lambda^{-1/2}X,X_n)=\frac{\chi_3}{(2k-2j)!}\Lambda^{m-k}(\nabla^{2k-2j}p_{m-j})(0,X_n)X^{2k-2j}+\\ \Lambda^{m-k-1/2}f_j(X,X_n,\Lambda)X^{2k-2j+1}.
\end{multline}
$f_j$ a son support contenu dans une région $W_3=\{(X,X_n) ; \Lambda^{-1/2}|X|\leq R\}$, répétant l'argument de (\ref{eq:3.5}), on voit que $f_j\in S(1,g_\rho)$, $X^{2k-2j+1}\in S(d_\rho^{2k-2j+1})$. Donc dans une région comme $X_1$, où $|X|\leq C_\rho$, on a un terme total dans $S(\Lambda^{m-k-1/2},g_\rho)$. $\rho$ est une constante indépendante de $\Lambda$, il faudra donc supposer que $\mu=\Lambda^{-1/2}\rho^{2k+1}$ est petit. Il n'y a pas de doute qu'on peut aussi supposer que $\chi_3$ vaut $1$ sur un voisinage du support de $\chi_1$.

Donc :
\begin{equation}
\chi_1^2\sum_{j=0}^k \Lambda^{m-j}p_{m-j}(\Lambda^{-1/2}X,X_n)=\chi_1^2\Lambda^{m-k}Q_k(X_n)(X)+S(\Lambda^{m-k}(\mu+\Lambda^{-1}),g_\rho).
\end{equation}
On fait maintenant trois constatations :
\begin{enumerate}[i)]
\item Pour chaque $X_n$, $Q_k$ a une borne inférieure positive sur $L^2$.
\item Pour chaque $X_n$,
\begin{equation} Q_k\geq cN_k-C\text{ dans } L^2 ;
\end{equation}
où $N_k=\sum_{|\alpha|\leq k}(L^{\alpha})^* L^{\alpha}$ où $L=(D_y+ i y)$ est un créateur de $L^2$. On sait que pour tout $s\in\R$, $N_k^s\in S(d_1^{2ks},g_1)$ et a pour symbole modulo $S(d_1^{2ks-2},g_1)$, la fonction $(\sum_j|L_j(x,\xi)|^2)^{ks}$. C'est par exemple une conséquence des travaux de J.M. Bony beaucoup plus généraux [2].
\item Les points i) et ii) entraînent qu'on a aussi $Q_k\geq cN_k$.
\end{enumerate}
On déduit de ces trois remarques qu'on peut appliquer à la $\chi_1^2Q_k$ une inégalité de sharp-Gärding vectorielle comme opérateur pseudo-différentiel en $X_n$ à valeurs dans $\LL(L^2)$. En effet pour chaque $X_n$, l'opérateur $N_k^{-1/2}Q_k(X_n)N_k^{-1/2}$ est défini positif et borné sur $L^2$ et il est aussi dans $S(1,g_1)$. L'opérateur $\chi_1(N_k^{-1/2}Q_k(X_n)N_k^{-1/2})\chi_1\geq c((\chi_1)^w)^2$, il y a donc $C$ telle que :
\begin{equation}\label{eq:3.15}
(\widetilde Q_k v,v)\geq -C\Lambda^{-1}\|v\|^2+c\|\chi_1 v\|^2.
\end{equation}
On applique (\ref{eq:3.15}) à $v=N_k^{1/2}u$, soit :
\begin{equation}(\widetilde Q_k N_k^{1/2}u,N_k^{1/2}u)\geq -C\Lambda^{-1}\|N_k^{1/2}v\|^2 +c\|\chi_1 N_k^{1/2} u\|^2.
\end{equation}
$N_k^{1/2}\widetilde Q_k N_k^{1/2}=(Q\chi_1^2)^{w_{1,\Lambda}}+R$, $R$ est un commutateur entre $N_k^{1/2}$ ou $Q_k$ et $\chi_1$, donc le calcul symbolique $S(m,g_1)\# S(m',g_\rho)$ s'applique (\ref{eq:3.2}), la fonction $h$ du calcul mixte vaut $h\leq d_\rho^{-1}d_1^{-1}$. On trouve que $r\in S(\rho^N d_\rho^{-N-1} d_1^{2k-1},g_1)$ pour tout $N\geq 0$, quand $N=0$ on obtient $r\in S(\rho^{-1} d_1^{2k-1},g_1)$. De même $[N_k^{1/2},\chi_1]\in S(d_1^{k-1} d_\rho^{-1},g_1)$. On a obtenu :
\begin{equation}\label{eq:3.17}
\Lambda^{-m+k}((Q\chi_1^2)^{w_{1,\Lambda}}\chi_0u ,\chi_0 u) \geq c\|N_k^{1/2}\chi_1\chi_0 u\|^2 -C\rho^{-1}\|N_{k-1}^{1/2}\chi_0 u\|^2-C\mu\|u\|^2.
\end{equation}
Il est évident qu'on peut remplacer $N^{1/2}$ par n'importe quel opérateur elliptique de $S(d_1^k,g_1)$.

Il suffit maintenant d'ajouter (\ref{eq:3.10}) et (\ref{eq:3.17}) pour obtenir :
\begin{equation}
\Lambda^{-m+k}((Q)^{w_{1,\Lambda}}\chi_0 u,\chi_0 u)\geq c\|d_1\chi_0 u\|^2 - C\rho^{-1}\|d_1^{k-2}\chi_0 u\|^2-C\mu\|u\|^2.
\end{equation}
On conclut en faisant tendre correctement $\rho\to \infty$ et $\mu\to 0$. On a donc prouvé le théorème.
\nocite{Bere,Ivri,Lasc,Lern,Lern2,Bony}


\textsc{Bernard Lascar.
 Richard Lascar. Université Denis Diderot. Département De Mathématiques.
Institut Mathématiques De Jussieu, Analyse Algébrique, 2 Place Jussieu, 75005 Paris. France}

{\it E-mail address :} {\tt richard.lascar@imj-prg.fr}
\end{document}